\documentclass[review]{amsart}

\usepackage{lineno,hyperref}
\usepackage[utf8]{inputenc}

\usepackage{hyperref}
\hypersetup{nesting=true,debug=true,naturalnames=true}
\usepackage{graphicx,amssymb,upref}

\usepackage{amsmath,amsthm}
\usepackage{amsfonts}
\usepackage{latexsym}
\usepackage{amsbsy}
\usepackage{amssymb,mathscinet}

\numberwithin{equation}{section}

\usepackage{amsmath,amsfonts,stackrel,amsthm,enumerate,subcaption}
\usepackage{graphicx,setspace,mathrsfs,placeins,mathtools,amssymb}
\usepackage{algorithm}
\usepackage{algorithmic,colortbl}
\newcommand{\algrule}[1][.2pt]{\par\vskip.5\baselineskip\hrule height #1\par\vskip.15\baselineskip}

\makeatletter
\newcommand\fs@norules{\def\@fs@cfont{\bfseries}\let\@fs@capt\floatc@ruled
	\def\@fs@pre{}%
	\def\@fs@post{}%
	\def\@fs@mid{\kern3pt}%
	\let\@fs@iftopcapt\iftrue}
\makeatother
\floatstyle{norules}
\restylefloat{algorithm}

\usepackage{amsfonts,amssymb,amscd,amsmath,enumerate,verbatim,calc,enumerate}
\theoremstyle{plain}
\newtheorem{theorem}{Theorem}[section]

\newtheorem{lemma}[theorem]{Lemma}
\newtheorem{proposition}[theorem]{Proposition}

\newtheorem{definition}[theorem]{Definition}

\newtheorem{remark}[theorem]{Remark}
\newtheorem{example}[theorem]{Example}

\modulolinenumbers[5]

\begin{document}

\title[L\'evy processes with lower incomplete gamma jumps]{L\'evy processes with jumps governed by lower incomplete gamma subordinator and its variations}

\author[]{Meena Sanjay Babulal$^*$$^a$, Sunil Kumar Gauttam$^*$$^b$ and Aditya Maheshwari$^\#$$^c$}
\email{$^a$19pmt002@lnmiit.ac.in, $^b$sgauttam@lnmiit.ac.in, $^c$adityam@iimidr.ac.in}

\address[]{
	$^*$Department of Mathematics, The LNM Institute of Information Technology, Rupa ki Nangal, Post-Sumel, Via-Jamdoli
	Jaipur 302031,
	Rajasthan, India. }
\address[]{$^\#$Operations Management and Quantitative Techniques Area,
	Indian Institute of Management Indore, Indore 453556, Madhya Pradesh, India.}
\begin{abstract}
In this paper, we study the L\'evy process 
time-changed by  independent L\'evy subordinators, namely, the 
incomplete gamma subordinator, the $\epsilon$-jumps incomplete 
gamma subordinator and tempered incomplete gamma 
subordinator.  We derive their important  distributional properties such as mean, variance, correlation, tail probabilities and fractional moments. The long-range dependence property of these processes are discussed. An application in insurance domain is studied in detail. Finally, we present the 
simulated sample paths for the subordinators.
\end{abstract}

 \maketitle

\section{Introduction}

	The stochastic models with random time clock appears in  various fields of applications such as finance (see \cite{Clark1973ASS,TCFPP-pub}), physics
(see \cite{Cahoy_2011,nezhadhaghighi2011first, STANISLAVSKY2008643,Zhang_2009}), ecology (see \cite{https://doi.org/10.1029/2001GL014123}), biology (see \cite{PhysRevLett.96.098102}) and \textit{etc}. As a result, there is ever increasing interest among probabilists into this kind of research problems, and it has  given arise to a new field of study called as the stochastic subordination.   Stochastic subordination involves investigating the stochastic process where the time variable is replaced by a non-decreasing L\'evy process. Its study can be divided  into two major classes, namely, the diffusion processes and the counting processes. 
A pioneer work on the stochastic subordination was first published by Bochner (see \cite{Bochner1949DiffusionEA, cooper_1957}) and subsequently many scholars studied various aspect of subordinated stochastic process such as homogeneity (see \cite{woll1959homogeneous}), Markov property (see \cite{1195520680}), long-range dependence (LRD) (see \cite{biardapp, maheshwari2016long}), and infinite divisibility( see \cite{applebaum2009levy, sato}), \textit{etc.} A comprehensive coverage can be found in  Bertoin (see \cite{bertoin1999subordinators}) and Sato (see \cite{sato1999levy}). In this paper, we focus on stochastic subordination of the general L\'evy  process.\\\\

Originally  Mandelbrot and Taylor (see \cite{Mandelbrot}) explored the idea of application of a subordinated process model in stock price . In their work, they examined a one-sided $\beta$-stable process with $0 < \beta < 1$ as the subordinator, effectively serving as a time-change mechanism for Brownian motion. They interpreted this subordinator as representing trading volume (or the number of transactions) up to time $t$. Later, several studies (see, for example, \cite{Barn-Niel97,Clark73, Heyde-Len05,Madan-Car-Chang98, Madan-Seneta90}) explored real-life use cases of the time-changed version of L\'evy processes. It shows wide interest and applicability in time-changed L\'evy processes. As a special case of the L\'evy process, the  Poisson  process is a well known and applicable model for count data. Buchak and Sakhno 
(see \cite{Buchak2017CompositionsOP}) investigated the Poisson process subordinated with gamma subordinator. Kumar \textit{et
	al.} (see \cite{Kumar-TCPP}) have discussed various characteristics of the Poisson process subordinated with the stable/inverse stable 
subordinator and the inverse Gaussian subordinator. Orsingher and Toaldo 
(see \cite{OrsToa-Berns}) explored the  Poisson 
process subordinated with a  L\'evy subordinator. 

Beghin and Ricciuti (see \cite{ricciuti}) defined the incomplete 
gamma (InG) subordinator, the incomplete gamma subordinator with jumps of size greater than or equal $\epsilon$ 
(InG-$\epsilon$) subordinator and tempered incomplete gamma (TInG) 
subordinator using lower-incomplete gamma function. The InG subordinator is defined as a non-decreasing L\'evy process with the Laplace exponent $\alpha\gamma(\alpha;\eta)$, where $\gamma(\alpha;\eta)$ is the lower-incomplete gamma function defined as

\begin{equation*}
	\gamma(\alpha;\eta)=\int_0^\eta e^{-y}y^{\alpha-1}dy,~~ \eta>0,0<\alpha\leq 1.
\end{equation*}The InG-$\epsilon$ subordinator is a modification of the InG subordinator whose jumps are greater than $\epsilon>0$ with the Laplace exponent $\tfrac{\alpha}{\epsilon^{\alpha}}\gamma(\alpha; \eta \epsilon)$. The TInG subordinator is defined as a non-decreasing L\'evy process with the Laplace exponent $\alpha \gamma(\alpha;\eta+\theta)-\alpha \gamma(\alpha;\theta) $, where $\theta>0$ is the tempering parameter. It has finite moment of any integer order. In this paper, we consider the InG, InG-$\epsilon$ and TInG subordinators as random clocks for the L\'evy  process. Our goal is to study important distributional properties, such as, mean, variance, correlation, tail probabilities and fractional moments.\\\\ 

The LRD property concerns with the memory of stochastic process. A stochastic model having the LRD or long ``memory" indicates that it is a non stationary process. This property can provide an alternative explanation to the empirical phenomenon that exhibits memory over a period of time; it has been investigated in detail 
 (see \cite{samorodnitsky2016stochastic}) and the references therein.  The definition of the LRD property is based on the second order property of stochastic processes; more specifically asymptotic behaviour of correlation function.  \\ \\ The classical approach in the insurance domain is to use the Poisson process to model the number of claim arrivals.  It is  noteworthy that the Poisson process does not possesses the LRD property. Moreover, only some time-changed Poisson processes exhibit the LRD property (\textit{eg. }the Poisson process time-changed by inverse stable or tempered stable subordinator (see \cite{biardapp,tsfnbp})); while others do not (\textit{e.g.} the Poisson process time-changed by gamma or stable subordinator). In our case, the InG and the InG-$\epsilon$ subordinators have infinite mean and therefore it is difficult to study a second order characteristic like LRD property. However, in the TInG subordinator case, we are able to find the results related  to the long-range behavior.  In this paper, we prove that the TInG subordinator and the L\'evy process subordinated with TInG subordinator have the LRD property. This is being used to study an application in ruin theory by taking Poisson process as a special case of the L\'evy process.\\\\ 
The Poisson process is used  to model risk for an insurance company. We have used the subordinated Poisson process (with the TInG  subordinator) as an alternative to the classical Poisson process in risk model for insurance, and is described as follows
\begin{equation*}
	Y(t)=ct-\sum_{j=1}^{N( S_{\alpha, \theta}(t))}X_{j}, t\geq 0,
\end{equation*}
where $c > 0$ denotes premium rate, which is assumed to be constant, $X_{j}$ be non-negative i.i.d. random variables with distribution $F$, representing the claim size and the Poisson process $\{N(t)\}_{t\geq0}$ subordinated by the TInG $\{ S_{\alpha, \theta}(t)\}_{t\geq0}$. We derive the governing equation for the joint probability that ruin happens in finite time and the deficit at the time of ruin. We also compute 
the joint distribution of ruin time and deficit at ruin when the initial capital is zero.\\\\ 
The simulation of sample paths provides a visual aid to understand a stochastic process.  We present the simulated sample paths for the InG, InG-$\epsilon$ and TInG subordinators. We have used the Metropolis algorithm(see \cite{casella2002statistical}) to simulate the sample paths where the candidate density is obtained by truncating the support of the exponential density; this approach is developed for the TInG subordinator. \\\\

The paper is organised as follows. In Section \ref{section 2}, we present some preliminary results that are required. Section \ref{section 3} deals with the asymptotic behaviour of tail probability  of L\'evy process subordinated by  InG and InG-$\epsilon$ subordinators. In Section \ref{section 4}, the asymptotic behaviour of fractional moment is presented for L\'evy process subordinated by InG, InG-$\epsilon$ and TInG subordinators. In Section \ref{section 5}, we derive the long range depedence property of TInG and L\'evy process subordinated by TInG subordinator. Section \ref{section 6} discusses  the  application of the subordinated Poisson process in insurance domain. In Section \ref{section 7}, we simulate the sample paths  of the InG,  InG-$\epsilon$ and TInG subordinators.

\section{Preliminaries}\label{section 2} 	\noindent In this section, we present some preliminary results which are required later in the paper. 

\setcounter{equation}{0}
Let $\mathbb{Z}^{+} :=\left\{0,1, 2,\cdots\right\}$ be the set of non-negative integers. Let $\left\{N (t, \lambda)\right\}_{t\geq0}$ be a Poisson process with rate $\lambda > 0$, so that
\begin{align*}
	p(n|t, \lambda):=\mathbb{P}[N (t, \lambda)=n]=\frac{(\lambda t)^{n}e^{-\lambda t}}{n!}, \qquad n \in \mathbb{Z}^{+}.
\end{align*}
For simplicity of notation we write $\left\{N (t, \lambda)\right\}_{t\geq0}$ as $\left\{N(t)\right\}_{t\geq0}$, when no confusion arises.\\
For $\alpha \in (0, 1]$, the  InG subordinator $\lbrace S_{\alpha}(t)\rbrace_{t \geq 0} $   (see \cite{ricciuti}) can be represented as a compound Poisson process
\begin{align*}
	S_{\alpha}(t) =
	\sum_{j=1}^{N_{\alpha}(t)}	Z_{j}^{\alpha},
\end{align*}
where $\left\lbrace N_{\alpha}(t)\right\rbrace_{ t \geq 0  }$
is a homogeneous Poisson process with the rate $\lambda := \alpha
\Gamma (\alpha) $ and the jumps $Z_{j}^{\alpha}$ are i.i.d. random variables, taking values in $[1, +\infty)$, with probability
density function
\begin{align}\label{Zpdf}
	f_{Z^{\alpha}}(z) = \frac{(z-1)^{-\alpha}z^{-1}1_{z\geq1}}{\Gamma (1-\alpha)\Gamma(\alpha)} = \frac{\sin(\pi \alpha)1_{z\geq1} }{\pi (z-1)^{\alpha}z },   \quad   \alpha \in (0,1).
\end{align}

When $\alpha = 1$, the jumps are unitary, and the subordinator $\lbrace S_{\alpha}(t)\rbrace_{t \geq 0} $  coincides with the Poisson process (see \cite{ricciuti}).	
Note that the subordinator $\lbrace S_{\alpha}(t)\rbrace_{t \geq 0} $ 
have jumps of size greater than or equal to $1$.\\
Similarly, the InG-$\epsilon$ subordinator $\lbrace S_{\alpha}^{(\epsilon)}(t)\rbrace_{t \geq 0}$ can be represented as a compound
Poisson process
\begin{equation*}
	S_{\alpha}^{(\epsilon)}(t) =\sum_{j=1}^{N^{\epsilon}(t)}	Z_{j}^{(\alpha,\epsilon)},
\end{equation*}
where $N^{\epsilon} := \lbrace N^{\epsilon}(t)\rbrace_{ t \geq 0 }$ 
is a homogeneous Poisson process with the rate $\lambda := \alpha
\Gamma (\alpha) \epsilon^{-\alpha} $ and the jumps $Z_{j}^{(\alpha,\epsilon)}$ are i.i.d. random variables, taking values in $[\epsilon, +\infty)$, with probability density function
\begin{align}\label{ZE pdf}
	f_{Z_{j}^{(\alpha,\epsilon)}}(z) = \frac{\epsilon^{\alpha}(z-\epsilon)^{-\alpha}z^{-1}1_{z\geq\epsilon}}{\Gamma (1-\alpha)\Gamma(\alpha)} 
	,   \quad   \alpha \in (0,1).
\end{align}
In contrast to the InG subordinator, the InG-$\epsilon$ subordinator $\lbrace S_{\alpha}^{(\epsilon)}(t)\rbrace_{t \geq 0}$ have jumps of size greater than or equal to $\epsilon$ (see \cite{ricciuti}). \\
The TInG subordinator $\lbrace S_{\alpha, \theta}(t)\rbrace_{t \geq 0}$ can be represented as a compound
Poisson process

\begin{equation*}
	S_{\alpha, \theta}(t)=	\sum_{j=1}^{N_{ \alpha, \theta}(t)}	Z_{j}^{ \alpha, \theta}, 
\end{equation*}
where $ N_{\alpha, \theta}$ := $\left\lbrace N_{\alpha, \theta}(t)\right\rbrace_{t\geq 0} $ is a homogeneous Poisson process with rate $\lambda := \alpha \Gamma(\alpha;\theta)$, where $\Gamma(\alpha;\theta)$ is the upper incomplete gamma function defined as
\begin{equation}\label{InGfU}
	\Gamma(\alpha;\eta)=\int_\eta^\infty e^{-y}y^{\alpha-1}dy,~~ \eta>0,0<\alpha\leq 1.\nonumber
\end{equation}

We have the following relationship between upper incomplete gamma function, lower incomplete gamma function and gamma function $$\Gamma(\alpha;\theta)+\gamma(\alpha;\theta)=\Gamma(\alpha).$$  The jumps $Z_{j}^{ \alpha, \theta}$
are i.i.d. random variables, taking values in $[1, +\infty)$ and with the probability density function
\begin{equation*}
	f_{Z_{j}^{ \alpha, \theta}}=\frac{e^{-\theta z}(z-1)^{-\alpha}z^{-1}1_{z\geq1}}{\Gamma (1-\alpha)\Gamma(\alpha;\theta)}, \quad \alpha \in (0,1).
\end{equation*} 
Observe that
the mean for InG and InG-$\epsilon$ subordinators does not exist, but  mean and variance of the TInG  subordinator $S_{\alpha, \theta}(t)$ are given by (see \cite{ricciuti})
\begin{equation*}
	\mathbb{E} S_{\alpha, \theta}(t) = t\alpha\theta^{\alpha-1}e^{-\theta}, \end{equation*}
\begin{equation*}
	\mbox{Var}  S_{\alpha, \theta}(t)= t\alpha\theta^{{\alpha-1}}e^{-\theta}+t({\alpha-1})\alpha\theta^{\alpha-2}e^{-\theta}.
\end{equation*}

The following result (see \cite{kumar2019fractional}) is key to our computation for the fractional order moments of the subordinators and subordinated L\'evy processes.
\begin{proposition}
    
	Let $X$ be a positive random variable with the Laplace transform $\widetilde{f(t)} $. Then its $q^{th}$ order moment, where $ q \in (n-1,n)$ is given by
	\begin{align}\label{laplace-erdelyi}
		\mathbb{E}(X^{q})=\frac{(-1)^n}{\Gamma{(n-q)}}\int_{0}^{\infty}\frac{d^n}{du^n}[\widetilde{f(u)}]u^{n-q-1} du.
	\end{align}  
\end{proposition}

\section{ asymptotic behaviour of tail probability } \label{section 3}
In this section, we study the asymptotic behaviour of tail probability of  L\'evy process subordinated with the InG and InG-$\epsilon$ subordinators. First we have the following definition.
\begin{definition}\label{def LInG}
	The L\'evy process subordinated with the InG subordinator (LInG) is defined as 
	\begin{align*}
		Q_\alpha(t) := Y( S_{\alpha}(t)), \qquad t \geq 0, 
	\end{align*}where $\lbrace {Y(t)}\rbrace_{t\geq 0} $ is the L\'evy process  and independent of the InG subordinator  $\left\{S_{\alpha}(t)\right\}_{t\geq0}$.
\end{definition}
If $\psi(\cdot)$ is the Laplace exponent of the L\'evy process $\{Y(t)\}_{t\geq0}$, then the Laplace exponent of the LInG 
$\{Q_\alpha(t)\}_{t\geq 0}$ is $e^{-t \alpha \gamma(\alpha;\psi(\eta))}.$  
 The asymptotic behaviour of the tail probability  of the LInG process is given by following theorem.
\begin{theorem}\label{tail prob of LevyInG}
	Let $\{Y(t)\}_{t\geq0}$ be a L\'evy process with the Laplace exponent $\psi(\cdot)$. 
 We assume that for  $\alpha \in(0,1)$ there exist  constants $ c\in \mathbb{R}$ and $\beta\in(0, 1)$ such that 
 \begin{align}\label{tail prob exponent eqn}
     \frac{\psi(\eta)^{\alpha}}{\eta} \sim c^{\alpha}\eta^{\beta-1},  \qquad \mbox{ as } \eta \rightarrow 0.
 \end{align}
 Then for any $t\geq 0$, we have 
	\begin{align}\label{levy asymping}
		\mathbb{P}(Q_\alpha(t)> x)\sim \frac{c^{\alpha} t  x^{-\beta}}{\Gamma(1-\beta)}, \mbox{ as $x\rightarrow \infty$ .}	\end{align}
	
\end{theorem}
\begin{proof} We consider the Laplace transform of tail probability of the $\{Q_\alpha(t)\}_{t\geq0}$,  
	for $\eta>0$
	\begin{align*}
		\int_{0}^{\infty}e^{-\eta x}\mathbb{P}(Q_\alpha(t)> x) dx
		&= \frac{1-\mathbb{E}e^{-\eta Y(S_{\alpha}(t))}}{\eta}\\
		&=\frac{1-e^{-t \alpha \gamma(\alpha;\psi(\eta))}}{\eta}.\end{align*}
	Using the Taylor's approximation up to first order, for $\eta\rightarrow 0$ we obtain
	\begin{align*}
		\frac{1-e^{-t \alpha \gamma(\alpha;\psi(\eta))}}{\eta}&\sim\frac{1-(1-t \alpha \gamma(\alpha;\psi(\eta)))}{\eta} \\
		&\sim\frac{t\alpha}{\eta}\frac{ \psi(\eta)^{\alpha}}{\alpha} \ \ \ \  \mbox{  $ \left( \gamma(\alpha;\eta)\sim \frac{\eta^{\alpha}}{\alpha}, \mbox{ as } \eta \rightarrow 0 \right) $   }\\
		&\sim t\frac{\psi(\eta)^{\alpha}}{\eta} \\
		&= c^{\alpha} \eta^{\beta-1} t.
	\end{align*}
	The desired result
	follows  from the Tauberian theorem (see  \cite[ p.446]{feller2008introduction})  for any $ t \geq 0$.\qedhere 
\end{proof} 
Now, we look at some important special cases of the above Theorem.
\begin{example}
      Let $\{N(t)\}_{t\geq 0}$ be a Poisson process with rate $\lambda>0$, 
    then for $\eta>0$ the Laplace exponent of $\{N(t)\}_{t\geq 0}$ is given by $\lambda (1-e^{-\eta})$.
     Taylor's approximation gives us  
    \begin{align*}
        \frac{(\lambda (1-e^{-\eta}))^{\alpha}}{\eta}\sim  \lambda^{\alpha}\eta^{\alpha-1}, \mbox{ as }\eta\to0.
    \end{align*}
    
     Comparing the above with $\eqref{tail prob exponent eqn}$, we get $c=\lambda^{}$ and $\beta=\alpha$. By Theorem \ref{tail prob of LevyInG}, asymptotic behaviour of  the tail probability of $\{N(S_{\alpha}(t))\}_{t\geq 0}$ is  
   \begin{align*}
		\mathbb{P}(N(S_{\alpha}(t))> x)\sim \frac{t \lambda^{\alpha} x^{-\alpha}}{\Gamma(1-\alpha)}, \mbox{  \qquad  as $x\rightarrow \infty$.} 	
  \end{align*}
\end{example}

\begin{example}
      Let $\{B(t)\}_{t\geq0}$ be a standard Brownian motion, then it has the Laplace exponent $\frac{\eta^{2}}{2}$ for $\eta>0$.
    We have 
    \begin{align*}
        \frac{(\frac{\eta^{2}}{2})^{\alpha}}{\eta}=  \frac{1}{2^{\alpha}}\eta^{2\alpha-1},\quad \forall \eta>0.
    \end{align*}
     Thus for $\alpha \in (0,\frac{1}{2})$, $c= \frac{1}{2^{}}$ and $\beta=2\alpha$ after comparing above with \eqref{tail prob exponent eqn}.  By Theorem \ref{tail prob of LevyInG}, for $\alpha \in (0,\frac{1}{2})$ the asymptotic behaviour of  the tail probability of $\{B(S_{\alpha}(t))\}_{t\geq 0}$ is  
   \begin{align*}
		\mathbb{P}(B(S_{\alpha}(t))> x)\sim \frac{t  x^{-2\alpha^{}}}{2^{\alpha}\Gamma(1-2\alpha^{})}, \mbox{\qquad as  $x\rightarrow \infty$.}
 \end{align*}
\end{example}

\begin{example}
      Let $\{D(t)\}_{t\geq 0}$ be a $\zeta$- stable subordinator with $\zeta\in (0,1)$, 
    then for $\eta>0$ the Laplace exponent of $\{D(t)\}_{t\geq 0}$ is given by $\eta^{\zeta}$.
    We have $$\frac{(\eta^{\zeta})^{\alpha}}{\eta}=  \eta^{\zeta\alpha-1},\quad \forall \eta>0.$$
     Identifying  $c=1$ and $\beta=\zeta\alpha$ with \eqref{tail prob exponent eqn}.  By Theorem \ref{tail prob of LevyInG}, asymptotic behaviour of  the tail probability of $\{D(S_{\alpha}(t))\}_{t\geq 0}$ is  
   \begin{align*}
		\mathbb{P}(D(S_{\alpha}(t))> x)\sim \frac{t  x^{-\zeta\alpha}}{\Gamma(1-\zeta\alpha)}, \mbox{\qquad  as $x\rightarrow \infty$.} 	
  \end{align*}
\end{example}

Now, we define  L\'evy process subordinated with the InG-$\epsilon$ subordinator. 

\begin{definition}\label{def LInG-e}
	The L\'evy process subordinated with the InG-$\epsilon$ subordinator (LInG-$\epsilon$) is defined as 
	\begin{align*}
		Q_{\alpha}^{(\epsilon)}(t) := Y( S_{\alpha}^{(\epsilon)}(t)), \qquad t \geq 0, 
	\end{align*}where $\lbrace {Y(t)}\rbrace_{t\geq 0} $ is the L\'evy process  and independent of the InG-$\epsilon$ subordinator  $\left\{S_{\alpha}^{(\epsilon)}(t)\right\}_{t\geq0}$.
\end{definition}
If $\psi(\cdot)$ is the Laplace exponent of the L\'evy processes  $\left\{Y(t)\right\}_{t\geq 0}$, then Laplace exponent of the L\'evy process $Y( S_{\alpha}^{(\epsilon)}(t))$ is $e^{-\frac{t \alpha}{\epsilon^{\alpha}} \gamma(\alpha;\psi(\eta))}.$
The  asymptotic behaviour of tail probability for the LInG-$\epsilon$ process is given by following theorem. The proof is similar to proof of Theorem \ref{tail prob of LevyInG}.
\begin{theorem}
	Let $\{Y(t)\}_{t\geq0}$ be a L\'evy process with the Laplace exponent $\psi(\cdot)$ satisfying \eqref{tail prob exponent eqn}.
 Then for any $t\geq 0$, we have 
	\begin{align*}
		\mathbb{P}(Q_{\alpha}^{(\epsilon)}(t))> x)\sim \frac{c^{\alpha} t  x^{-\beta}}{\Gamma(1-\beta)}, \mbox{\qquad  $x\rightarrow \infty$.} 	\end{align*}
	
\end{theorem}

\section{ Asymptotic behaviour of fractional moments }\label{section 4}
In this section, we investigate the asymptotic behaviour of fractional moments of L\'evy process subordinated by InG, InG-$\epsilon$ and TInG subordinators. First, we look at  the asymptotic behaviour of the fractional moments of LInG.

\begin{theorem}\label{gen-frac-moment-of-levy}
Let $\{Y(t)\}_{t\geq0}$ be a L\'evy process with the Laplace exponent $\psi(\cdot)$ with the series expansion $\psi(\eta)=\sum_{k=1}^{\infty}d_{k}\eta^{k}$. 
Let $\xi$ is the first non-zero coefficient $d_{k}$ and $m=k$ 
corresponding to first non-zero $d_{k}$ in expansion of $\psi(\eta)$ respectively. Then, for $p \in (0,1)$ the fractional moment of $p^{th}$ order of the process the $\left\{Q_\alpha(t)\right\}_{t\geq0}$  is finite for $p<\alpha$ and its asymptotic behaviour is  given by 
\begin{align*}
    \mathbb{E}[Q_\alpha(t)]^p\sim 
	\frac{\alpha c_{0} t^{\frac{p}{m\alpha} }\Gamma\left(1-\frac{p}{m\alpha}\right)}{\Gamma(1-p)}, \mbox{ as } t\to\infty,
\end{align*}
	where $c_{0}$ depends on $\alpha$, $\xi$, $m$ and $p.$
\end{theorem}
\begin{proof}
	We first argue the existence of fractional moments. The asymptotic behaviour of the tail probability given by \eqref{levy asymping} allows us to conclude
	$\mathbb{E}[Y( S_{\alpha}(t))]^p < \infty $ for $p < \alpha$.
 \noindent For the asymptotic behaviour of the fractional moments, 
	 using (\ref{laplace-erdelyi}) we get
	\begin{align*}
		\mathbb{E}[Q_\alpha(t)]^p&=\frac{-1}{\Gamma(1-p)}\int_{0}^{\infty}\frac{d}{d\eta}[e^{-t \alpha \gamma(\alpha;\psi(\eta))}]\eta^{-p}d\eta \nonumber \\ 
		&=\frac{1}{\Gamma(1-p)}\int_{0}^{\infty}e^{-t \alpha \gamma(\alpha;\psi(\eta))} \dfrac{d}{d\eta}\left[\alpha t  \gamma(\alpha;\psi(\eta))\right]\eta^{-p}d\eta \nonumber \\ 
  &=\frac{\alpha t}{\Gamma(1-p)}\int_{0}^{\infty}e^{-t \alpha \gamma(\alpha;\psi(\eta))} e^{-\psi(\eta)}\psi(\eta)^{\alpha -1}\psi'(\eta)\eta^{-p}d\eta. \nonumber 
	\end{align*} 
 \noindent Here 
	$ h(\eta)=\alpha \gamma(\alpha;\psi(\eta))$ and $\phi(\eta)=e^{-\psi(\eta)}\psi(\eta)^{\alpha -1}\psi'(\eta)\eta^{-p}$. Then, we have  
\begin{align*}
h(\eta) &=\alpha \gamma(\alpha;0) +\alpha\sum_{j=0}^{\infty}\frac{( \psi(\eta))^{\alpha+j}}{(\alpha+j)j!}, \mbox{ and }\\
 \phi(\eta)&=e^{-\psi(\eta)}\psi(\eta)^{\alpha -1}\psi'(\eta)\eta^{-p}\\
  &=\sum_{j=0}^{\infty}\frac{\psi(\eta)^j}{j!}\psi(\eta)^{\alpha -1}\psi'(\eta)\eta^{-p}\\
& = \sum_{j=0}^{\infty}\frac{\psi(\eta)^{j+\alpha-1}}{j!}\psi'(\eta)\eta^{-p},
	\end{align*}
	where $h(0)=\alpha \gamma(\alpha;0)$, $a_{0}=\xi_{}^{\alpha}$,   $\mu=m \alpha $, $b_{0}=\xi_{}^{\alpha} m $ and $\rho=m \alpha-p$. Now, we apply the Laplace–Erdelyi Theorem (see \cite{wojdylo2006coefficients}) to the above integral and we get
	\begin{align*}
		\mathbb{E}[Q_\alpha(t)]^p\sim \frac{\alpha t}{\Gamma(1-p)}\sum_{i=0}^{\infty}\frac{c_{i}}{t^{\frac{\rho+i}{\mu}}}\Gamma\left(\frac{\rho+i}{\mu}\right).
	\end{align*}
	Above series is dominated by first term for large $t$, 
	which leads to
	\begin{align*}
		\mathbb{E}[Q_\alpha(t)]^{p}\sim  \alpha c_{0} {t^{1-\frac{\rho}{\mu}}}\frac{\Gamma\left(\frac{\rho}{\mu}\right)}{\Gamma(1-p)}
\sim\frac{\alpha c_{0} t^{\frac{p}{m \alpha} }\Gamma\left(1-\frac{p}{m \alpha}\right)}{\Gamma(1-p)} ,
	\end{align*}
	where $c_{0} =\frac{b_{0}}{m \alpha a_{0}^{1-\frac{p}{m \alpha}}}=\frac{\xi^{\frac{p}{m}}}{\alpha}$. \qedhere
	\end{proof}

\noindent Now, we look at particular cases of the above theorem.


\begin{example}
      Let $\{N(t)\}_{t\geq0}$ be a Poisson process with rate $\lambda$, 
    then for $\eta>0$ the Laplace exponent of $\{N(t)\}_{t\geq 0}$ is given by $\psi(\eta)=\lambda (1-e^{-\eta})=\lambda(\eta-\frac{\eta^2}{2!}+\frac{\eta^3}{3!}+\cdots)$.
  Here first nonzero coefficient $\xi$ equal to $\lambda$  and  $m=1$. Then by Theorem \ref{gen-frac-moment-of-levy}, for $p \in (0,1)$ 
	 asymptotic behaviour of fractional moment of PInG is  given by 
   \begin{align*}\label{PInG asymping}
	\mathbb{E}[N( S_{\alpha}(t))]^p\sim \frac{(\lambda t)^{\frac{p}{\alpha}} \Gamma(1-\frac{p}{\alpha})}  {{\Gamma(1-p)}},\qquad \mbox{as  $t\rightarrow \infty$.}
  \end{align*}
\end{example}

\begin{example}
      Let $\{B(t)\}_{t\geq0}$ be a standard Brownian motion with the Laplace exponent $\psi(\eta)=\frac{\eta^{2}}{2}$ for $\eta>0$.
    Here we have first nonzero coefficient $\xi$ equal to $\frac{1}{2}$ and  $m=2$. Then by Theorem \ref{gen-frac-moment-of-levy}, for $p \in (0,\frac{1}{2})$     
	 asymptotic behaviour of  fractional moment of subordinated Brownian motion $\left\{B(S_{\alpha}(t))\right\}_{t\geq0}$ is given by 
   \begin{align*}
		\mathbb{E}[B( S_{\alpha}(t))]^p\sim \frac{  \Gamma(1-\frac{p}{2\alpha})}{2^{\frac{p}{2\alpha}}\Gamma(1-p)}t^{\frac{p}{2\alpha}},\mbox{\qquad as  $t\rightarrow \infty$}.
 \end{align*}
\end{example}
Using the same line arguments as in the proof of Theorem \ref{gen-frac-moment-of-levy}, we have the following result for LInG-$\epsilon$.
\begin{theorem} 
Let $\{Y(t)\}_{t\geq0}$ be a L\'evy process with the Laplace exponent $\psi(\cdot)$,
where  $\psi(\eta)=\sum_{k=1}^{\infty}d_{k}\eta^{k}$.Let $\xi$ is the first non-zero coefficient $d_{k}$ and $m=k$ 
corresponding to first non-zero $d_{k}$ in expansion of $\psi(\eta)$ respectively.
Then, for $p \in (0,1)$ the fractional moment of $p^{th}$ order of the process the $\left\{Y(S_{\alpha}^{\epsilon}(t))\right\}_{t\geq0}$  exists, is finite for $p<\alpha$ and 
	its asymptotic behaviour is  given by 
\begin{align*}
    \mathbb{E}[Q_{\alpha}^{(\epsilon)}(t)]^p\sim 
	\frac{\alpha c_{0} t^{\frac{p}{m\alpha} }\Gamma\left(1-\frac{p}{m\alpha}\right)}{\Gamma(1-p)},
\end{align*} 
	where $c_{0} =\frac{b_{0}}{\left(m \alpha {a_{0}^{1-\frac{p}{m \alpha}}}\right)}=\frac{\xi^{\frac{p}{m}}}{\alpha}$. \qedhere
\end{theorem}

Now, 
we
define the L\'evy process subordinated with the TInG subordinator and study asymptotic behaviour of fractional moments of TInG and   L\'evy process subordinated with the TInG subordinator.
\begin{definition}\label{def TInG}
	The L\'evy process subordinated with the TInG  subordinator (LTInG) is defined as 
	\begin{align*}
		Q_{\alpha,\theta}(t) = Y( S_{\alpha,\theta}(t)), \qquad t \geq 0, 
	\end{align*}
	where $\lbrace{Y(t)}\rbrace_{t\geq0} $ is the L\'evy process with rate $\lambda>0$ and independent of the TInG subordinator $\left\{S_{\alpha,\theta}(t)\right\}_{t\geq0}$.  
\end{definition}
If $\psi(\cdot)$ is the Laplace exponent of the L\'evy processes $\lbrace{Y(t)}\rbrace_{t\geq0} $, then the Laplace exponent of the L\'evy process $\lbrace{Y(S_{\alpha,\theta}(t))}\rbrace_{t\geq0} $ is $\exp({-t \alpha (\gamma(\alpha;\psi(\eta)+\theta)-\gamma(\alpha,\theta)}).$

Now, we state the results regarding  asymptotic behaviour of  fractional moments for the TInG subordinator and LTInG, which can be obtained along the similar line as in proof of Theorem  \ref{gen-frac-moment-of-levy}. 
\begin{theorem}
	
  Let $p \in (0,1]$, then  asymptotic behaviour of fractional moment of $p$-th order of the process  TInG  exist, is finite for $p < \alpha$ and is given by
		\begin{align*}
			\mathbb{E}[S_{\alpha,\theta}(t)]^{p}\sim  (\alpha e^{-\theta}\theta^{\alpha-1}t)^{p}, \qquad \mbox{ as } t\rightarrow \infty.
		\end{align*}
\end{theorem}  
\begin{theorem}\label{gen frac moment of levy TInG}
Let $\lbrace{ Y(t)}\rbrace_{t\geq0} $ be a L\'evy process with the Laplace exponent $\psi(\cdot)$,
where  $\psi(\eta)=\sum_{k=1}^{\infty}d_{k}\eta^{k}$. 
Then, for $p \in (0,1)$ the fractional moment of $p^{th}$ order of the LTInG $\left\{Q_{\alpha,\theta}(t)\right\}_{t\geq0}$ is finite for $p<\alpha$ and 
	its asymptotic behaviour is  given by 
\begin{align*}
    \mathbb{E}[Q_{\alpha,\theta}(t)]^p\sim 
	\frac{ (\alpha \theta^{\alpha-1} e^{-\theta} \xi t)^{\frac{p}{m} }\Gamma\left(1-\frac{p}{m}\right)}{\Gamma(1-p)},\qquad \mbox{ as } t\rightarrow \infty,
\end{align*} 
	 where $\xi$ and $m$ are the first non-zero coefficient $d_{k}$ and power of that term in expansion of $\psi(\eta)$ respectively.
\end{theorem}

\noindent Now we look at some important special cases of the above Theorem.

\begin{example}
    Let $\lbrace{ N(t)}\rbrace_{t\geq0} $ be a Poisson process with rate $\lambda$ and the Laplace exponent $\psi(\eta)=\lambda (1-e^{-\eta})$ for $\eta>0$.
 Then  by Theorem \ref{gen frac moment of levy TInG}, for  $p \in (0,1)$ 
	 asymptotic behaviour of fractional moment of PTInG  is  given by 
   \begin{align*}
	\mathbb{E}[N( S_{\alpha,\theta}(t))]^p\sim 
	{ (\alpha \theta^{\alpha-1} e^{-\theta} \lambda t)^{{p} }},\qquad \mbox{ as } t\rightarrow \infty.
  \end{align*}
\end{example}
\begin{example}
      Let $\lbrace{ B(t)}\rbrace_{t\geq0} $ be a standard Brownian motion with the Laplace exponent $\psi(\eta)=\frac{\eta^{2}}{2}$ for $\eta>0$.
   Then  by Theorem \ref{gen frac moment of levy TInG}, for $p \in (0,1)$ 
	 asymptotic behaviour of fractional moment of subordinated Brownian TInG 
  is  given by

   \begin{align*}
    \mathbb{E}[B( S_{\alpha,\theta}(t))]^p\sim 
	\frac{ (\alpha \theta^{\alpha-1} e^{-\theta}  t)^{\frac{p}{2} }\Gamma\left(1-\frac{p}{2}\right)}{2^\frac{p}{2}\Gamma(1-p)},\qquad \mbox{ as }t\rightarrow \infty.
\end{align*} 
\end{example}

\section{long range depedence }\label{section 5}
In this section, we discuss long range depedence property of the the TInG  subordinator and L\'evy process subordinated by TInG subordinator.
We first state the definition of LRD  (see \cite{maheshwari2016long}).
\begin{definition}\label{lrd defn}
	Let $0 < s < t$ and $s$ be fixed. Assume a stochastic process $\left\{X(t)\right\}_{t\geq0}$
	has the correlation function $\mbox{Corr}\left[X(s),X(t)\right]$ that satisfies
	$$c_{1}(s)t^{-d}\leq \mbox{Corr}\left[X(s),X(t)\right] \leq c_{2}(s)t^{-d}$$
	for large $t, d > 0$, $c_{1}(s) > $0 and $c_{2}(s) > 0$. That is, 
	
	$$\lim_{t\rightarrow \infty}\frac{\mbox{Corr}\left[X(s),X(t)\right]}{t^{-d}}=c(s)$$
	for some $c(s) > 0$ and $d > 0$. We say $\left\{X(t)\right\}_{t\geq0}$ has the long-range dependence
	(LRD) property if $d \in \left(0,1\right)$ and has the short-range dependence (SRD) property if
	$d \in \left(1,2\right)$.
\end{definition}

Now we show that the TInG  $\left\{S_{\alpha, \theta}(t)\right\}_{t\geq 0}$ has LRD property.
\begin{theorem}
	The TInG $\left\{S_{\alpha, \theta}(t)\right\}_{t\geq 0}$ has LRD property.
\end{theorem}

\begin{proof}
	First we compute the covariance using independent increment property of subordinator.	For $0\leq s<t<\infty$, we have
	\begin{align}
		\mbox{Cov}[S_{\alpha,\theta}(s),S_{\alpha,\theta}(t)]&= \mbox{Cov}[S_{\alpha, \theta}(s),(S_{\alpha,\theta}(t)-S_{\alpha,\theta}(s))+S_{\alpha,\theta}(s)]\nonumber\\
		&=\mbox{Cov}[S_{\alpha,\theta}(s),(S_{\alpha,\theta}(t)- S_{\alpha, \theta}(s))]+ \mbox{Cov}[S_{\alpha,\theta}(s),S_{\alpha,\theta}(s)]\nonumber\\
		&=\mbox{Var}[S_{\alpha,\theta}(s)].
  \end{align}
	Thus the correlation function is given by 
	\begin{align}
		\mbox{Corr}[S_{\alpha, \theta}(s),S_{\alpha, \theta}(t)]&=\frac{\mbox{Cov}[S_{\alpha, \theta}(s),S_{\alpha, \theta}(t)]}{\mbox{Var}[S_{\alpha, \theta}(s)]^{1/2}\mbox{Var}[S_{\alpha, \theta}(t)]^{1/2}}\nonumber\\
		&=\frac{\mbox{Var}[S_{\alpha, \theta}(s)]^{1/2}}{\mbox{Var}[S_{\alpha, \theta}(t)]^{1/2}}\nonumber\\
		&=s^{1/2}t^{-1/2}. \label{corr of tempered subo}
	\end{align}
	Hence
	$$\lim_{t\rightarrow\infty}\frac{\mbox{Corr}[ S_{\alpha, \theta}(s), S_{\alpha, \theta}(t)]}{t^{\frac{-1}{2}}}=s^{\frac{1}{2}}.$$\label{LRD of tempered subor}
	Therefore, the TInG $\left\{S_{\alpha, \theta}(t)\right\}_{t\geq 0}$ has LRD property.	\end{proof}

Since the TInG has finite mean and variance, it is expected that the LTInG also have finite mean and
variance. Following is a special case of a more general result (see \cite[Theorem 2.1]{Leonenko2014CorrelationSO} ).

 \begin{lemma}\label{LTInG m, var and cov}
     Let $\lbrace{Y(t)}\rbrace_{t\geq0} $  be a L\'evy process with finite second moment. Then the mean, variance and covariance of the $\lbrace{Q_{\alpha,\theta}(t)}\rbrace_{t\geq0} $  are given by 
	\begin{itemize}
		\item[(a)] $\mathbb{E}[Q_{\alpha,\theta}(t)]= t \alpha\theta^{\alpha-1}e^{-\theta} \mathbb{E}[Y(1)]$,
		\item[(b)] $\mbox{Var}[Q_{\alpha,\theta}(t)]=\mathbb{E}[Y(1)]^2(t\alpha\theta^{{\alpha-1}}e^{-\theta}+t({\alpha-1})\alpha\theta^{\alpha-2}e^{-\theta})+\mbox{Var}[Y(1)](t\alpha\theta^{\alpha-1}e^{-\theta})$,
		\item[(c)] Cov$[Q_{\alpha,\theta}(s),Q_{\alpha,\theta}(t)]=\mbox{Var}[Y(1)](s \alpha\theta^{\alpha-1}e^{-\theta} )+\mathbb{E}[Y(1)]^{2}(s\alpha\theta^{{\alpha-1}}e^{-\theta}+s({\alpha-1})\alpha\theta^{\alpha-2}e^{-\theta})$, for $0\leq s<t$.
	\end{itemize}
 \end{lemma}

\begin{theorem}\label{LTING-lrd}
	The LTInG $\left\{Q_{\alpha,\theta}(t)\right\}_{t\geq 0}$ has the LRD property.
\end{theorem}
\begin{proof}
	Using the covariance of  the LTInG  from Lemma \ref{LTInG m, var and cov}, we derive expression for correlation function of the LTInG   as
	\begin{align}
		&\mbox{Corr}[Y( S_{\alpha, \theta}(s)),Y( S_{\alpha, \theta}(t))]\nonumber\\
		&=\frac{\mbox{Cov}[Y( S_{\alpha, \theta}(s)),Y( S_{\alpha, \theta}(t))]}{(\mbox{Var}[Y( S_{\alpha, \theta}(s))])^{\frac{1}{2}}(\mbox{Var}[Y( S_{\alpha, \theta}(t))])^{\frac{1}{2}}}\nonumber\\
		&=t^{\frac{-1}{2}}s^{\frac{1}{2}}\nonumber.
	\end{align}
	Hence
	\begin{align*}
		\lim_{t\rightarrow\infty}\frac{\mbox{Corr}[Y( S_{\alpha, \theta}(s)),Y( S_{\alpha, \theta}(t))]}{t^{\frac{-1}{2}}}=s^{\frac{1}{2}}.
	\end{align*}
	This completes the proof of LRD property for $\left\{Y(S_{\alpha, \theta}(t))\right\}_{t\geq 0}$.\qedhere
 \end{proof}

\noindent Now we look at some important special cases.

\begin{example}\label{PTInG m, var and cov}
Let $\lbrace{N(t)}\rbrace_{t\geq0} $ be a Poisson process with rate $\lambda>0, \mathbb{E}[N(1)]=\lambda$ and $\mbox{Var}[N(1)]=\lambda$. Then the mean, variance and covariance of $\lbrace{ N( S_{\alpha, \theta}(t))}\rbrace_{t\geq0} $ are given by 
	\begin{itemize}
		\item[(a)] $\mathbb{E}[N( S_{\alpha, \theta}(t))]=\lambda t\alpha\theta^{\alpha-1}e^{-\theta}$,
		\item[(b)] $\mbox{Var}[N( S_{\alpha, \theta}(t))]=\lambda^2(t\alpha\theta^{{\alpha-1}}e^{-\theta}+t({\alpha-1})\alpha\theta^{\alpha-2}e^{-\theta})+\lambda(t\alpha\theta^{\alpha-1}e^{-\theta})$,
		\item[(c)] Cov$[N( S_{\alpha, \theta}(s)), N( S_{\alpha, \theta}(t))]=\mbox{Var}[N( S_{\alpha, \theta}(s))]$, for $0\leq s<t$.
	\end{itemize}
 By Theorem $\ref{LTING-lrd}$, the PTInG $\left\{N(S_{\alpha, \theta}(t))\right\}_{t\geq 0}$ has the LRD property.
 \ifx
The Laplace exponent of L\'evy process  $\lbrace{ N( S_{\alpha, \theta}(t))}\rbrace_{t\geq0} $ is  $\exp\left({-t \alpha [\gamma(\alpha;\lambda (1-e^{\eta})+\theta)-\gamma(\alpha;\theta)]}\right),$
and the pgf is given by $\exp\left({-t\alpha(\gamma(\alpha;\lambda(1-u)+\theta)-\gamma(\alpha;\theta))}\right).$
The pmf for $\lbrace{ N( S_{\alpha, \theta}(t))}\rbrace_{t\geq0} $ is given by
\begin{align*}
	\mathbb{P}(N( S_{\alpha, \theta}(t))=0)&=\exp{\left(-t\alpha\rho(\alpha;\theta)\right)},\\
	\mathbb{P}(N( S_{\alpha, \theta}(t))=1)&=((\lambda+\theta)^{\alpha-1}\alpha\lambda t ) \exp{\left(-t\alpha\rho(\alpha;\theta)\right)}, \\
	\mathbb{P}(N( S_{\alpha, \theta}(t))=2)&=\frac{1}{2}
	\left[{\left(1-\alpha\right)} \alpha {\left(\theta + \lambda\right)}^{\alpha - 2} \lambda^{2} t + {\left(\alpha {\left(\theta + \lambda\right)}^{\alpha - 1} \lambda t \exp{\left(-\theta - \lambda\right)} + \lambda\right)} \alpha {\left(\theta + \lambda\right)}^{\alpha - 1} \lambda t \right] e^{ -\alpha t {\rho(\alpha;\theta)} - \theta - \lambda},\\
	\mathbb{P}(N( S_{\alpha, \theta}(t))=3)&=\frac{1}{6} \left[(\alpha - 1) {(\alpha - 2)} \alpha {(\theta + \lambda)}^{\alpha - 3} \lambda^{3} t - 2 \, {\left(\alpha {\left(\theta + \lambda\right)}^{\alpha - 1} \lambda t e^{-\theta - \lambda} + \lambda\right)} {\left(\alpha - 1\right)} \alpha {\left(\theta + \lambda\right)}^{\alpha - 2} \lambda^{2} t 
	\right.\nonumber \\ & \left.\qquad +{\left(\alpha {\left(\theta + \lambda\right)}^{\alpha - 1} \lambda t e^{-\theta - \lambda} + \lambda\right)}^{2} \alpha {\left(\theta + \lambda\right)}^{\alpha - 1} \lambda t
	-
	\left({\left(\alpha - 1\right)} \alpha {\left(\theta + \lambda\right)}^{\alpha - 2} \lambda^{2} t e^{-\theta - \lambda} \right.\right.\nonumber\\ & \ \ \ \ \ \ \ - \left. \left. \alpha {\left(\theta + \lambda\right)}^{\alpha - 1} \lambda^{2} t e^{-\theta - \lambda}\right)\alpha {\left(\theta + \lambda\right)}^{\alpha - 1} \lambda 
	\right] \exp{\left(-\alpha t {\rho(\alpha;\theta)} - \theta - \lambda\right)},
\end{align*}
where $\rho(\alpha;\theta)=\gamma\left(\alpha; \theta + \lambda\right) - \gamma\left(\alpha; \theta\right)$ and so on.
The pmf of the the PTInG 
	is also given by 
	\begin{align}
		\mathbb{P}[N( S_{\alpha, \theta}(t))=n] =\frac{\lambda^{n}}{n!}\mathbb{E}[e^{-\lambda  S_{\alpha, \theta}(t)}( S_{\alpha, \theta}(t))^{n}],\quad  n=0,1,2,3\cdots. \label{PTInG pmf form 2}
	\end{align}
	
\ifx 
\begin{proof}
	Let $g(y,t)$ be the probability density function of $S_{\alpha, \theta}(t)$. Then
	\begin{align*}
		\mathbb{P}[N( S_{\alpha, \theta}(t))=n]&=\int_{0}^{\infty}\mathbb{P}[N( S_{\alpha, \theta}(t))=n\big|S_{\alpha, \theta}(t)]g(y,t)dt\\
		&=\int_{0}^{\infty}\frac{(\lambda y)^{n}e^{-\lambda y}}{n!}g(y,t)dt\\
		&=\frac{\lambda^{n}}{n!}\mathbb{E}[e^{-\lambda  S_{\alpha, \theta}(t)}( S_{\alpha, \theta}(t))^{n}]. \qedhere
	\end{align*}	
	
\end{proof}
\fi
The representation (\ref{PTInG pmf form 2}) allows easy verification of  the normalizing condition
$$\sum_{n=0}^{\infty}\mathbb{P}[N( S_{\alpha, \theta}(t))=n]=1.$$
Consider \begin{align*}
	\sum_{n=0}^{\infty}\mathbb{P}[N( S_{\alpha, \theta}(t))=n]&=\sum_{n=0}^{\infty}\frac{\lambda^{n}}{n!}\mathbb{E}[e^{-\lambda  S_{\alpha, \theta}(t)}( S_{\alpha, \theta}(t))^{n}]\\
	&=\sum_{n=0}^{\infty}\int_{0}^{\infty}\frac{(\lambda y)^{n}e^{-\lambda y}}{n!}g(y,t)dy\\
	&=\int_{0}^{\infty}g(y,t)dy=1.\qedhere
\end{align*}
Using simple algebraic calculations, one can see that the transition probabilities of
the $\lbrace{ N( S_{\alpha, \theta}(t))}\rbrace_{t\geq0} $ are given by
\begin{align} &\mathbb{P}[N( S_{\alpha, \theta}(t+h))=n|N( S_{\alpha, \theta}(t))=m] \nonumber \\&=
	\left\{ \begin{array}{ll}\label{transition probability of tpp}
		1-\alpha f(\lambda)
		+o(h) & \mbox{if $n=m$  } \\  
		-h\left[(-1)^{i}\frac{\lambda^{i}}{i!}f^{(i)}(\lambda)\right]+o(h)  & \mbox{if $ n=m+i, i=1,2,3 \cdots$, }\\
	\end{array}\right.
\end{align}
where $f(\lambda)=	\gamma(\alpha;\lambda+\theta)-\gamma(\alpha;\theta)$ is Laplace exponent of $S_{\alpha, \theta}(t)$	.

\fi

\end{example}
\begin{example}\label{BM m, var and cov}
	Let $\lbrace{ B(t)}\rbrace_{t\geq0} $ be a standard Brownian motion, then $\mathbb{E}[B(1)]=0$ and $\mbox{Var}[B(1)]=1$. The mean, variance and covariance of $\lbrace{ B( S_{\alpha, \theta}(t))}\rbrace_{t\geq0} $  are given by 
	\begin{itemize}
		\item[(a)] $\mathbb{E}[B( S_{\alpha, \theta}(t))]=0$,
		\item[(b)] $\mbox{Var}[B( S_{\alpha, \theta}(t))]=t\alpha\theta^{\alpha-1}e^{-\theta}$,
		\item[(c)] Cov$[B( S_{\alpha, \theta}(s)),B( S_{\alpha, \theta}(t))]=\mbox{Var}[B( S_{\alpha, \theta}(s))]$, for $0\leq s<t$.
	\end{itemize}
 By Theorem $\ref{LTING-lrd}$, the BTInG $\left\{B(S_{\alpha, \theta}(t))\right\}_{t\geq 0}$ has the LRD property.
\end{example}

\section{Application in insurance ruin}\label{section 6}
  The primary focus of life insurance companies revolves around the management of mortality and longevity risk. Long-memory processes are very useful statistical models that excel at capturing persistent dependencies and correlations across extended time spans. In insurance domain, the existence of long memory is established (see \cite{yan_peters_chan_2021}) in the mortality data or death count data. The authors show that forecasts from models without a long memory structure provide overestimated mortality rates which will give rise to underestimated life expectancy. This work also provides a new approach to enhance mortality forecasts in terms of accuracy and reliability.  In \cite{WangSAJ,WANG202125}, authors demonstrate application of non-Markovian model with LRD to help insurers to mitigate risks associated with longevity or mortality in the life insurance market. They prove that the LRD
has a significant effect on longevity hedging, and suggest reinsurance as an actuarial risk management tool is robust to the LRD property of the mortality rate.\\\\
The ruin theory is a branch of actuarial science that deals with the financial modeling of the likelihood of a company or individual becoming insolvent. The classical risk process of insurance defined below models the  distribution of claims,  balance of assets and liabilities over time 
$$Z(t)= ct-\sum_{j=0}^{N(t)}X_{j}, t\geq 0,$$
where  $c > 0$ is fixed premium rate and $\left\{N(t)\right\}_{t\geq 0}$ is the homogeneous Poisson process	which counts claims arrival till time $t$. The claim amount $X_{j}$ with distribution $F$ is independent of $\left\{N(t)\right\}_{t\geq 0}$.

We here propose to use the PTInG process $\{N(S_{\alpha,\theta}(t)\}_{t\geq 0}$ replacing the Poisson process $\{N(t)\}_{t\geq 0}$ in the classical risk process $\{Z(t)\}_{t\geq0}$. The number of claims in $\{Z(t)\}_{t\geq0}$ follows the Poisson distribution which assumes that the arrivals are i.i.d. while our proposed model has the LRD property (see Example \ref{PTInG m, var and cov}).

The LRD property assumes some sort of dependence on the successive claims and it is a more closer approximation of a real-life situation.   \\\\
Consider the risk model 
\begin{equation}\label{ risk model}
	Y(t)=ct-\sum_{j=1}^{N( S_{\alpha, \theta}(t))}X_{j}, t\geq 0,
\end{equation}
where $c > 0$  denotes a constant premium rate and $X_{j}$ are non-negative i.i.d. random variables with distribution $F$, representing the claim size.

The joint probability of ruin and deficit is a measure used in actuarial science to assess the financial stability of an insurance company. It describes the probability that an insurance company will not only become insolvent, or ``ruined," but also that it will have a deficit in its reserves. 
This measure is used to evaluate the effectiveness of different risk management strategies, such as adjusting pricing, increasing reserves, or purchasing reinsurance. Actuaries use this measure to evaluate the overall financial stability of the company and to make decisions on how to manage its risks.

\noindent In this section, we derive results for the ruin probability, joint distribution of time to ruin and deficit at ruin, and derive its  governing differential equation for our proposed model  \eqref{ risk model}.  
Note that the transition probabilities of the $\{N(S_{\alpha,\theta}(t)\}_{t\geq 0}$ are given by
\begin{align} &\mathbb{P}[N( S_{\alpha, \theta}(t+h))=n|N( S_{\alpha, \theta}(t))=m] \nonumber \\&=
	\left\{ \begin{array}{ll}\label{transition probability of tpp}
		1- hf(\lambda)
		+o(h) & \mbox{if $n=m$  } \\  
		-h\left[(-1)^{i}\frac{\lambda^{i}}{i!}f^{(i)}(\lambda)\right]+o(h)  & \mbox{if $ n=m+i, i=1,2,3 \cdots$, }\\
	\end{array}\right.
\end{align}
where $f(\lambda)=	\gamma(\alpha;\lambda+\theta)-\gamma(\alpha;\theta)$ is the Laplace exponent of $\left\{S_{\alpha, \theta}(t)\right\}_{t\geq 0}$.

The premium loading factor, denoted by $\rho$, signifies the profit margin of the insurance firm and is defined as the following ratio 
$$\rho=\frac{\mathbb{E}[Y(t)]}{\mathbb{E}\left[\sum_{j=1}^{N( S_{\alpha, \theta}(t))}X_{j} \right]}=\frac{ct}{\mu \mathbb{E}[N( S_{\alpha, \theta}(t))]}-1,$$
where $\mu= \mathbb{E}[X_{j}].$  Let us denote the initial capital by $u > 0$. Define the surplus process $\{U(t)\}_{t\geq 0}$ by 
$$U(t)=u+Y(t), ~~~t\geq 0.$$
The insurance company will be called in ruin if the surplus process falls below the zero
level. Let $T_u$ be the random variable which denotes the first time to ruin. It is defined as
$$T_u = \inf\left\{t>0 :U(t)<0\right\}.$$
The probability of ruin is given by
$\psi(u)=\mathbb{P}\left\{T_u<\infty\right\}.$
The joint probability that ruin happens in finite time and the deficit at the time of
ruin, denoted as $D=|U(T_{u})|$, is given by
\begin{equation}\label{ruin prob}
	G(u,y)=\mathbb{P}\left\{T_u<\infty, D\leq y\right\},\quad y\geq 0.
\end{equation}

Using (\ref{transition probability of tpp}), we get 
\begin{align*}
	G(u,y)&=	(1- h f(\lambda))G(u+ch,y)\\
	&~~~~~~-h\frac{(-1)\lambda f'(\lambda)}{1!}\left[\int_{0}^{u+ch}G(u+ch-x,y)dF(x)+(F(u+ch+y)-F(u+ch))\right]\\
	&~~~~~~-h\frac{(-1)^{2}\lambda^{2} f''(\lambda)}{2!}\left[\int_{0}^{u+ch}G(u+ch-x,y)dF(x)+(F(u+ch+y)-F(u+ch))\right]\\
	&~~~~~~~-\cdots\\
	&=(1- h f(\lambda))G(u+ch,y)-h\sum_{n=1}^{\infty}\frac{(-\lambda)^{n}}{n!}f^{(n)}(\lambda)\times\\&~~~~~~~\left[\int_{0}^{u+ch}G(u+ch-x,y)dF(x)+(F(u+ch+y)-F(u+ch))\right]
\end{align*}

After rearranging the terms, we have that 
\begin{align*}
	\frac{G(u+ch,y)-G(u,y)}{ch}&=\frac{1}{c}f(\lambda)G(u+ch,y)+\left(\frac{1}{c}\sum_{n=1}^{\infty}\frac{(-\lambda)^{n}}{n!}f^{(n)}(\lambda)\left[\int_{0}^{u+ch}G(u+ch-x,y)dF(x)\right.\right.\\& \qquad+\left .\left(F(u+ch+y)-F(u+ch)\right)\bigg]\right).\end{align*}

Now taking limit $h\rightarrow 0$, we get
\begin{align*}
	\frac{\partial G}{\partial u}&=\frac{f(\lambda)}{c}G(u,y)+\left(\frac{1}{c}\sum_{n=1}^{\infty}\frac{(-\lambda)^{n}}{n!}f^{(n)}(\lambda)\left[\int_{0}^{u}G(u-x,y)dF(x)+(F(u+y)-F(u))\right]\right)\\
	&=\frac{f(\lambda)}{c}G(u,y)+\left(\frac{1}{c}\sum_{n=1}^{\infty}\frac{(-\lambda)^{n}}{n!}f^{(n)}(\lambda)\left[\int_{0}^{u}G(u-x,y)dF(x)+(F(u+y)-F(u))\right]\right),
\end{align*}

Using Taylor's series, we get
\begin{align*}
	\sum_{n=1}^{\infty}\frac{(-\lambda)^{n}}{n!}f^{(n)}(\lambda)&=\sum_{n=0}^{\infty}\frac{(-\lambda)^{n}}{n!}f^{(n)}(\lambda)-f(\lambda)\\
	&=f(0)-f(\lambda)\\
	&=-f(\lambda).
\end{align*}
Thus, we have proved the following result.
\begin{theorem}
	Let $G(u,y)$ defined in (\ref{ruin prob}), denote the joint probability distribution of
	time to ruin and deficit at this time of the risk model (\ref{ risk model}). Then, it satisfies the following
	integro-differential equation
	\begin{align}\label{join prob ruin thm}
		\frac{\partial G}{\partial u}&=\frac{f(\lambda)}{c}\left[G(u,y)-\int_{0}^{u}G(u-x,y)dF(x)-(F(u+y)-F(u))\right].
	\end{align}
\end{theorem}
The rate of change in $G(u, y)$ with respect to change in initial capital $u$ can be found using the integro-differential equation \eqref{join prob ruin thm}.
 We observe from \eqref{join prob ruin thm}
that  $G(u,y)$ has  inverse relationship with initial capital $u$.\\

The solution of integro-differential equation $G(u,y)$,  defines corresponding density $g(u, y)$ as $g(u, y)=\frac{d}{dy}G(u,y)$ (see \cite{article3}). Thus $g(u, y) dy$ is the probability that ruin occurs and $Y(t)$ will be between $-y$ and $-y+dy$. The adjustment coefficient appear as one of the root when we solve an equation to compute the Laplace transform of $g(u,y)$ with respect to `$u$' (see \cite{article3}). The adjustment coefficient is closely related to ruin probability (see \cite[Theorem 13.4.1]{Bowers1997ActuarialMS}).\\

Above equation \eqref{join prob ruin thm} is difficult to solve as it contains both derivative and integral of $G(u,y)$. So to get rid of derivative of $G(u,y)$, we integrate it to get following result.

\begin{theorem}
	The joint distribution of ruin time and deficit at ruin when the initial capital is zero, $G(0,y)$ is given by
	\begin{align}\label{join prob ruin thm at zero initial capital }
		G(0,y)&=\frac{f(\lambda)}{c}\left[\int_{0}^{\infty}(F(u+y)-F(u))du\right].
	\end{align}
\end{theorem}
\begin{proof}
	On integrating (\ref{join prob ruin thm}) with respect to $u$ on $(0, \infty)$, we obtain
	\begin{align*}
		G(\infty,y)-G(0,y)&=\frac{f(\lambda)}{c}\left[\int_{0}^{\infty}G(u,y) du -\int_{0}^{\infty}\int_{0}^{u}G(u-x,y)dF(x)du  \right.\nonumber \\
  &\left.\ \ \ \ \ \ \ \ \ \ \ \  \ \ \ \ \ -\int_{0}^{\infty}(F(u+y)-F(u)du)\right].
		\shortintertext{Note that $G(\infty, y) = 0$, then}
		G(0,y)&=\frac{f(\lambda)}{c}\left[\int_{0}^{\infty}(F(u+y)-F(u))du\right].\qedhere
	\end{align*}
\end{proof}

\begin{remark}
In $G(u,y)$, $y$ denotes the limited liability, when ruin will happen. By taking $y\to \infty$ in $G(u,y)$, we shift liability from  finite to infinite. As $\lim_{y\to\infty}G(u,y)=\psi(u)$,  where $\psi(u)$ denotes the ruin probability with infinite liability.
	 As $y\to \infty$ in (\ref{join prob ruin thm at zero initial capital }), we get
	\begin{equation*}
		\psi(0)=\frac{f(\lambda)}{c}\left[\int_{0}^{\infty}(1-F(u))du\right]=\frac{\mu f(\lambda)}{c}.	\end{equation*}

This shows that ruin probability $\psi(u)$ at $u=0$ depends only on the expected claim size $\mu$ and not on specific form of claim size distribution $F(u)$. On taking $y\to \infty$ in (\ref{join prob ruin thm}), we obtain
	\begin{align*}
		\frac{\partial \psi}{\partial u}&=\frac{f(\lambda)}{c}\left[\psi(u)-\int_{0}^{u}\psi(u-x)dF(x)-(1-F(u))\right].
	\end{align*}

\end{remark}
 Similar to the integro-differential equation \eqref{join prob ruin thm}, we obseve that ruin probability $\psi(u)$  also has inverse relationship with initial capital.
\section{Simulation }\label{section 7}
In this section, we present  algorithms to simulate sample paths for the
InG, the InG-$\epsilon$, the TInG subordinators.

To simulate the InG subordinator, we first
calculate the cdf $F_{Z^{\alpha}}(\cdot)$  
of random variable  $Z_{}^{\alpha} $ with pdf given by \eqref{Zpdf}. Let $F_{Z^{\alpha}}(\cdot)$ be the cdf of  random variable  $Z_{}^{\alpha} $, then  
\begin{align}\label{cdf of Zj}
	F_{Z^{\alpha}}(x)=\int_{-\infty}^{x}	f_{Z^{\alpha}}(z)dz &= \int_{-\infty}^{z}\frac{(z-1)^{-\alpha}z^{-1}1_{z\geq1}}{\Gamma (1-\alpha)\Gamma(\alpha)}dz= \int_{-\infty}^{x}\frac{\sin(\pi \alpha)1_{z\geq1} }{\pi ((z-1)^{-\alpha})z }dz \nonumber  \\
	&=\int_{1}^{x}\frac{(z-1)^{-\alpha}z^{-1}}{\Gamma (1-\alpha)\Gamma(\alpha)}dz=\frac{\pi \csc(\pi \alpha)-I_{\frac{1}{x}}(\alpha, (1-\alpha))}{\Gamma(\alpha)\Gamma(1-\alpha)}\nonumber\\
	&={1-B_{\frac{1}{x}}(\alpha, (1-\alpha))},
\end{align}

where $B_{x}(a,b)=\int_{0}^{x} y^{a-1} (1-y)^{b-1}dy$ is the incomplete beta function.

Now,  we present algorithm to simulate the InG 
subordinator
using algorithm of compound Poisson process  (see \cite{abdullah2020algorithm}).
\algrule\vspace*{-.4cm}
\begin{algorithm}[H]
	\caption{Simulation of the InG subordinator}\label{algo-ing}
	\begin{algorithmic}[1]
		\vspace*{-.22cm}
		\algrule
		\renewcommand{\algorithmicrequire}{\textbf{Input:}}
		\renewcommand{\algorithmicensure}{\textbf{Output:}}
		\REQUIRE $\lambda>0 $,
		$\alpha\in(0,1)$ and $T\geq0$.
		\ENSURE  $Y(t)$,  simulated sample paths for the InG subordinator.
		\\ \textit{Initialisation} : $t=0$ and $Y=0$.
		\WHILE{$t < T$} 
		\STATE generate a uniform random variable $U \sim U(0, 1)$.
		\STATE set $t \leftarrow t - U/\lambda$.
		\STATE generate i.i.d. random variable $Z^{\alpha}$ using the inverse transform method to cdf \eqref{cdf of Zj}.
		\STATE set $Y \leftarrow Y + Z^{\alpha}$.
		\ENDWHILE
		
		\RETURN $Y.$
		\algrule
	\end{algorithmic}\label{simu of CPP InG}
\end{algorithm}

\noindent The cdf $F_{Z^{\alpha,\epsilon}}$ of random variable $Z_{}^{(\alpha,\epsilon)}$  can be obtained on same line as \eqref{cdf of Zj} and by using \eqref{ZE pdf}
\begin{align}\label{cdf ofZe}
	F_{Z_{}^{(\alpha,\epsilon)}}(x)= {1-B_{\frac{\epsilon}{x}}(\alpha, (1-\alpha))}.
\end{align}
We can simulate  the sample paths for InG-$\epsilon$ subordinator using Algorithm \ref{algo-ing}  by replacing equation \eqref{cdf of Zj}   
by  \eqref{cdf ofZe} in Step 4:.

It can be noted that we can not  simulate the TInG  subordinator using Algorithm \ref{simu of CPP InG}  as the cdf of random variable $Z^{\alpha,\theta}$  
does not have closed form 
and therefore we turn to 
Metropolis 
algorithm which is a 
special case of the MCMC scheme. In this method, we generate  $Y$ with pdf $f_{Y}$, called as target density, by choosing another random variable $V$ with pdf $f_{V}$, called as candidate density, such that $f_{Y}$ and $f_{V}$ have common support. The process is repeated for a large number of iterations, and the resulting sequence of accepted values approximates the desired probability distribution.

We present the Metropolis algorithm (see \cite[p.$254$]{casella2002statistical}) with $f_{Z^{\alpha, \theta}}$ as the target density.  
We define the candidate density $f_{V}$ as
\begin{align}\label{candidate density tempered}
	f_{V} (v)=
	\left\{\begin{array}{ll}
		\frac{\lambda e^{-\lambda v}}{e^{-\lambda}}, &  v\in[1,\infty) \\
		0,   &  \mbox{otherwise}
	\end{array},\right.
\end{align}
by truncating the exponential density. To generate $Z^{\alpha, \theta} \sim f_{Z^{\alpha, \theta}}$ using $ f_{V} (v)$ we use Metropolis algorithm. Let $Z^{\alpha,\theta} \sim f_{Z^{\alpha, \theta}}(z)$ and $V \sim f_{V} (v)$, where $f_{Z^{\alpha, \theta}}$ and $f_{V}$ have
common support. 


\algrule\vspace*{-.3cm}

\begin{algorithm}[H]
	 
	\caption{Metropolis algorithm}\label{Metropolis}
	\begin{algorithmic}[1]
		\vspace*{-.22cm}
		\algrule
	
  \renewcommand{\algorithmicrequire}
  {\textbf{Input:}}
		\renewcommand{\algorithmicensure}{\textbf{Output:}}
		\REQUIRE  $f_{Z_{\alpha,\theta}}(z)$ and $f_V(v)$ with common support, parameter $\alpha$, $\theta$, and number of iterations $N$.
		\ENSURE  random variable with pdf $f_{Z_{\alpha,\theta}}(z)$.
		\STATE generate an initial sample $Z_0$ from the distribution $f_V(v)$.
		
		\FOR {$i = 1$ to N}
		\STATE generate a random variable $U_i \sim \text{uniform}(0,1)$ and $V_i \sim f_V$.
		\STATE calculate the acceptance probability $\rho_i \leftarrow \min\left\{ \frac{f_{Z_{\alpha,\theta}}(V_i)}{f_V(V_i)} \cdot \frac{f_V(Z_{i-1})}{f_{Z_{\alpha,\theta}}(Z_{i-1})}, 1 \right\}$.
		\STATE set \begin{equation*}
			Z_{i}	\leftarrow	\left\{ \begin{array}{ll}
				V_{i} & \mbox{if $U_{i} \leq \rho_{i}$  } \\  
				Z_{i-1}	  & \mbox{if $U_{i}> \rho_{i}$ }
			\end{array},\right.
		\end{equation*}
		\ENDFOR
		\RETURN $Z_{N}$
		\algrule
	\end{algorithmic}
\end{algorithm}
Above algorithm  produces random variables $Z_{i}$ from the pdf (approximately) $f_{Z^{\alpha,\theta}}$. 
We use this $Z_{i}$ random variables to simulate the TInG  subordinator.
\algrule\vspace*{-.4cm}
\begin{algorithm}[H]

 \caption{Simulation of the TInG  subordinator}\label{algo-ting}
	\begin{algorithmic}[1]
		\vspace*{-.22cm}
		\algrule
		\renewcommand{\algorithmicrequire}{\textbf{Input:}}
		\renewcommand{\algorithmicensure}{\textbf{Output:}}
		
		\REQUIRE $\lambda>0 $, $\alpha\in(0,1)$ and $T\geq 0$.\\
		\ENSURE  sample paths of $Y(t)$, the TInG subordinator.
		\\ \textit{Initialisation} : $t=0, Y=0$.
		\WHILE {$t \le T$ }
		\STATE generate a uniform random variable $U \sim U(0, 1)$.
		\STATE set $t \leftarrow t - (1/\lambda) * (U)$.
		\STATE generate i.i.d. random variable $Z$ using Algorithm \ref{Metropolis}.
		\STATE set $Y \leftarrow Y + Z$.
		\ENDWHILE
		\RETURN $Y$.
		\algrule
	\end{algorithmic}
\end{algorithm}
\begin{remark}  
The sample paths of the InG and the InG-$\epsilon$ subordinators can also be generated  using the Metropolis algorithm. To simulate the InG subordinator  we use the candidate density  \eqref{candidate density tempered}, and for  the InG-$\epsilon$ subordinator, the candidate density is defined as 
\begin{align*}
	f_{V} (v)=
	\left\{\begin{array}{ll}
		\frac{\lambda e^{-\lambda v}}{e^{-\lambda \epsilon}}, &  v\in[\epsilon,\infty) \\
		0,   &  \mbox{otherwise}
	\end{array}.\right.
\end{align*}
\end{remark}
Now, we look at the Blumenthal-Getoor index to investigate the jump distribution of sample path. We first state the definition of the Blumenthal-Getoor index (see \cite{todorov2011volatility}).

\begin{definition}
    Let $\left\{X(t)\right\}_{t\geq0}$ be a L\'evy process with L\'evy measure $\nu$, then  the Blumenthal-Getoor index $\beta$ of $\left\{X(t)\right\}_{t\geq0}$ is defined as  
\begin{equation*}
    \beta=\inf \left\{ p>0 : \int_{|x|\leq1}|x|^{p}\nu\left(dx\right)<\infty\right\}.
\end{equation*}
\end{definition}

The InG subordinator has the L\'evy triplet $(0,0,\pi)$ (see \cite{ricciuti}), where density of the L\'evy measure $\pi(dz)$ is $ \frac{\alpha(z -1)^{-\alpha}z^{-1}1_{z\geq1}}{\Gamma (1-\alpha)}.$
Hence 

\begin{align*}
  \int_{|z|\leq1}|z|^{p}\pi\left(dz\right)&= \int_{|z|\leq1}|z|^{p}\frac{\alpha(z-1)^{-\alpha}z^{-1}1_{z\geq1}}{\Gamma (1-\alpha)}dz= \int_{\{z=1\}}|z|^{p}\frac{\alpha(z-1)^{-\alpha}z^{-1}}{\Gamma (1-\alpha)}dz=0,
\end{align*}
for all $p>0.$ Therefore the Blumenthal-Getoor index $\beta$ for InG subordinator is $0.$\\

 Similarly, the support of L\'evy density corresponding to TInG subordinator is $[1,\infty)$. Hence the Blumenthal-Getoor index $\beta$ for TInG subordinator is also $0.$ \\
 
Note that the InG-$\epsilon$ subordinator has the L\'evy triplet $(0,0,\pi_{\epsilon})$ (see \cite{ricciuti}), where density of the L\'evy measure $\pi_{\epsilon}(dz)$ is  $ \frac{\alpha(z-\epsilon)^{-\alpha}z^{-1}1_{z\geq\epsilon}}{\Gamma (1-\alpha)}.$
Therefore 
\begin{align*}
  \int_{|z|\leq1}|z|^{p}\pi_{\epsilon}\left\{dz\right\}&= \int_{|z|\leq1}|z|^{p}\frac{\alpha(z-\epsilon)^{-\alpha}z^{-1}1_{z\geq\epsilon}}{\Gamma (1-\alpha)}dz\nonumber\\
  &=\frac{\alpha}{\Gamma (1-\alpha)}\int_{\epsilon}^{1}z^{p-1}(z-\epsilon)^{-\alpha}dz<\infty,
\end{align*}
for all $p>0.$ Hence the Blumenthal-Getoor index $\beta$ for InG-$\epsilon$ subordinator is $0.$

\ifx 
The TInG subordinator has the L\'evy triplet $(0,0,\pi_{\theta})$, where $\pi_{\theta}$ is the L\'evy measure $ \frac{\alpha(z-1)^{-\alpha}z^{-1}e^{-\theta z}1_{z\geq1}}{\Gamma (1-\alpha)}.$
Consider 
\begin{align}
  \int_{|z|\leq1}|z|^{p}\pi_{\theta}\left\{dz\right\}&= \int_{|z|\leq1}|z|^{p}\frac{\alpha(z-1)^{-\alpha}z^{-1}e^{-\theta z}1_{z\geq1}}{\Gamma (1-\alpha)}dz=0, 
\end{align}
for all $p>0.$ Therefore the Blumenthal-Getoor index $\beta$ for TInG subordinator is $0.$
\fi

 Since L\'evy measure of InG, InG-$\epsilon$ and TInG subordinators is finite(see \cite{ricciuti}), they all have finite jump activity (see \cite{Matsuda2004IntroductionTM}).
\subsubsection*{\bf Interpretation of sample paths} 
\begin{figure}[h]
	\begin{subfigure}{4cm}
		\centering
		\includegraphics[scale=.2]{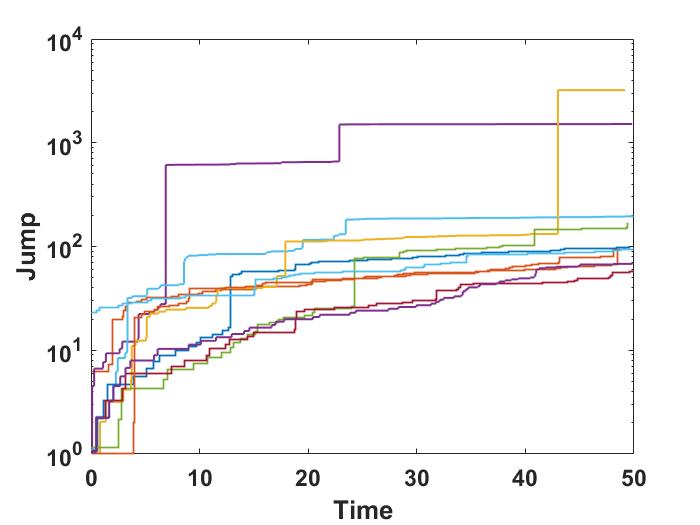}
		\subcaption{ Parameter $\alpha=0.8$ and $\lambda=\alpha*\Gamma(\alpha)$
			$=0.9314$}\end{subfigure}
	\hfill
\begin{subfigure}{4cm}
		\centering
	
	\includegraphics[scale=.2]{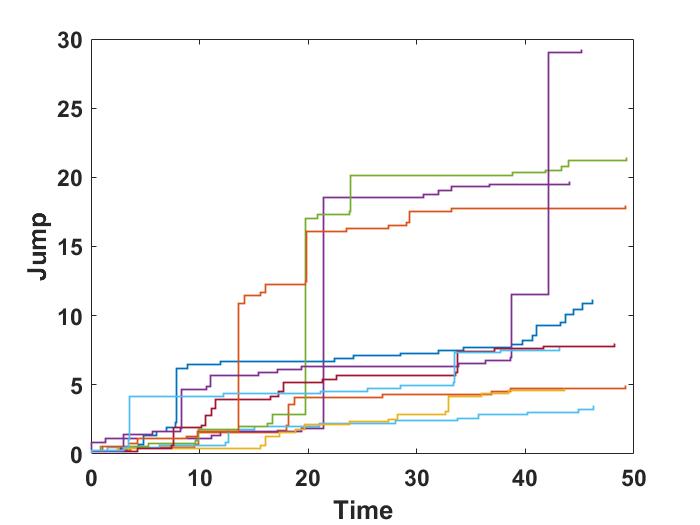}
		\subcaption{ Parameter : $ \alpha=0.8, \epsilon=0.2$ and $\lambda=0.2570$
   }
		
		
	\end{subfigure}
 \hfill
 \begin{subfigure}{4cm}
		\centering
		\includegraphics[scale=.2]{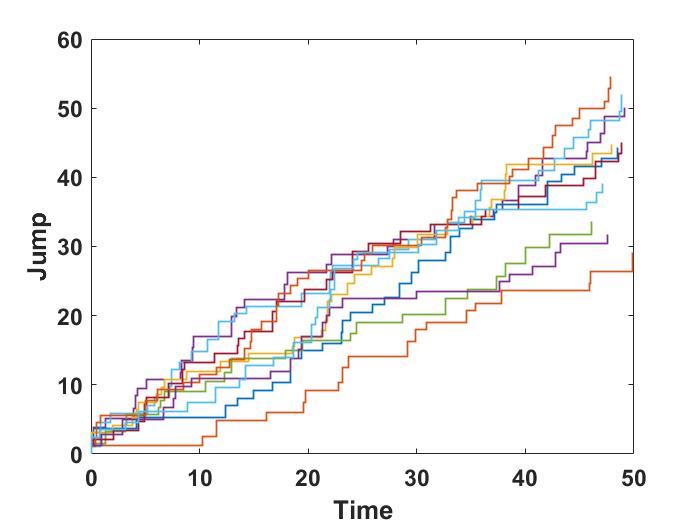}
		\subcaption{Parameter: $ \alpha=0.2$, $\theta=0.1$ and $\lambda=
  0.5802$
  }
		
	\end{subfigure}
 \caption{Ten simulated sample paths of the InG, InG-$\epsilon$ and TInG subordinators for different parameters.  }\label{fig subordinator}
\end{figure}

We simulate the samples paths of the InG and the InG-$\epsilon$ subordinators in Figures \ref{fig subordinator}(A) and \ref{fig subordinator}(B) respectively, using Algorithm \ref{algo-ing}. The sample paths of the TInG subordinator are presented in Figure \ref{fig subordinator}(C) using the Algorithm \ref{algo-ting}.
The InG subordinator has all the jumps greater than one while we can choose the jumps greater equal to $\epsilon$ in the InG-$\epsilon$ subordinator. It can be observed that, in comparison with other subordinators, the jump activity of the TInG subordinator is quite muted due to the tempering parameter $\theta$.
Our sample path not only reflects the finite jump activity of subordinators but also differentiates among their jump activity. It helps us to choose a model for jump processes. The parameter estimation of these subordinators will be an important and interesting problem to consider for a future work in this direction.\\

\paragraph{$\mathbf{Acknowledgment}$} We thank anonymous reviewer for carefully examining the paper which led to significant improvement in quality. First author would like to acknowledge the Centre for Mathematical \& Financial Computing 
and the DST-FIST grant for the infrastructure support for the computing 
lab facility under the scheme FIST (File No: SR/FST/MS-I/2018/24) 
at the LNMIIT, Jaipur.

\def\cprime{$'$}

\end{document}